\documentclass[12pt,reqno]{amsart}

\usepackage{color}
\usepackage{amsmath, amsthm, amssymb}
\usepackage{amsfonts}
\usepackage[ansinew]{inputenc}
\usepackage[dvips]{epsfig}
\usepackage{graphicx}
\usepackage[english]{babel}
\usepackage{hyperref}

\usepackage{cite}
\usepackage{graphicx}
\usepackage{amscd}
\usepackage{color}
\usepackage{bm}
\usepackage{enumerate}

\usepackage{verbatim}
\usepackage{hyperref}
\usepackage{amstext}
\usepackage{latexsym}
\theoremstyle{plain}
\newtheorem{Theorem}{Theorem}[section]

\newtheorem{Remark}[Theorem]{Remark}

\numberwithin{Theorem}{section} \numberwithin{equation}{section}

\def\square{\vbox{
\hrule height .4pt \hbox{\vrule width .4pt height 7pt \kern 7pt
\vrule width .4pt} \hrule height .4pt }}
\def\QED{\hfill {$\square$}\goodbreak \medskip}

\newcommand{\average}{{\mathchoice {\kern1ex\vcenter{\hrule height.4pt
width 6pt depth0pt} \kern-9.7pt} {\kern1ex\vcenter{\hrule height.4pt
width 4.3pt depth0pt} \kern-7pt} {} {} }}

\def\R{\mathbb{R}}

\def\div{\text{div}}

\newcommand{\e }{\varepsilon }

\renewcommand{\l }{\lambda }

\newcommand{\n }{\nabla }

\renewcommand{\phi}{\varphi}

\newcommand{\s }{\sigma }

\renewcommand{\o }{\omega }
\renewcommand{\O }{\Omega }

\newcommand{\ov}{\overline}

\newcommand{\be}{\begin{equation}}
\newcommand{\ee}{\end{equation}}

\newcommand{\de}{\partial}

\renewcommand{\H}{{\mathcal H}}

\newcommand{\cM}{{\mathcal M}}

\newcommand{\M}{\cM}

\renewcommand{\phi}{\varphi}

\renewcommand{\epsilon}{\varepsilon}

\begin{document}
\title{Existence of self-Cheeger sets on Riemannian Manifolds}
\author{Ignace Aristide Minlend}
\address{\ African Institute for Mathematical Sciences (A.I.M.S.) of Senegal KM 2, Route de Joal, B.P. 1418 Mbour, S\'en\'egal }
\email{ignace.a.minlend@aims-senegal.org}
\thanks{}

\begin{abstract}
Let  $(\M,g)$  be a smooth compact Riemannian manifold of dimension
$N\geq 2$. We prove the  existence of a family $(\O_\e)_{\e\in
(0,\e_0)}$ of  self-Cheeger sets in  $(\M,g)$. The domains
$\O_\e\subset\M$ are perturbations of geodesic balls of radius $\e$
centered at $p \in \M$, and in particular, if $p_0$ is a
non-degenerate critical point of the scalar curvature  of $g$, then
the family $(\de \O_\e)_{\e\in (0,\e_0)}$ constitutes a smooth
foliation of a neighborhood of $p_0$.
\end{abstract}
\maketitle

\section{Introduction and main results}
Given  a measurable set $A$ and an open set $\O$ in a compact Riemannian manifold
$(\M,g)$, the  relative perimeter, $P(A, \O)$, of $A$ in $\O$ is defined (see for
instance \cite{MR}) by
$$P(A, \O):=\sup\left\{\int_{A} \div_g(\xi)\,\textrm{d}v_g\,:\, \xi\in
C_c^1(\O),\, |\xi|\leq 1\right\}, $$ where $\textrm{d}v_g$ is the
Riemannian volume element and $C_c^1(\O)$ is the set of $C^1$ vectors fields with
 compact support in $\O$. The perimeter, $P(A)$, of a measurable set $A \subset \M$
 is the relative perimeter   $P(A, \M)$ of $A$ in $\M$.  If a set $A$ is of finite perimeter, denoting by
$\partial^{*}A$ the reduced boundary of $A$, $P(A)$  and $P(A, \O)$
are respectively the $(N-1)$-dimensional Riemannian volume of
$\partial^{*}A$ and $\partial^{*}A\cap \O$.

Let $u\in L^1(\O)$, the total variation of $u$ in $\O$ is defined by
$$|D u|(\O)=\sup\left\{\int_{\O} u \div_g(\xi)\,\textrm{d}v_g\,:\, \xi\in
C_c^1(\O),\, |\xi|\leq 1\right\}$$ and we say that $u$ is  a
function with bounded variation in $\O$ if its total variation is
finite. The space of functions with bounded variation in $\O$ is
denoted $BV(\O)$. We notice that $P(A, \O)$ is the total variation
of the characteristic function $1_A$ in $\O$.

 Let  $A$  be a set of
finite perimeter in $\M$, then  from De Giorgi's structure Theorem
\cite[Theorem 2.2]{A.FigalliandGe} (see also \cite[Remark on
p.161]{Hofmann} or \cite{Alter}), we have the Gauss-Green formula
\begin{equation}\label{The div}
\int_{A} \div_g(\xi)\,\textrm{d}v_g=\int_{\partial^{*}A} \langle \xi,\nu_A \rangle \textrm{d}\s_g\quad \textrm{for all} \quad\xi\in C^{1}(\M),
\end{equation}
where $\nu_A$ is the unit outer normal to $\de^{*} A$ and
$\textrm{d}\s_g$ is the area  element induced by the metric $g$.

Let $\O$
be an open set  in $(\M,g)$  with smooth boundary. The Cheeger
constant $h(\O)$  of $\O$ is defined by
\begin{equation}\label{eq: Cheeger}
h(\O):=\inf_A \frac{P(A)}{|A|},
\end{equation}
 where the domain $A$
varies over all measurable subsets of $\O$ with finite perimeter and
positive $N$-dimensional Riemannian volume  denoted  by $|A|$. It
was first introduced by J. Cheeger in \cite{JeffChe} to obtain
geometric lower bound of eigenvalues on $\O$.  We refer the reader
to \cite{Alter} and \cite{E Parini} for an introductory survey on
the Cheeger problem and to \cite{BB,Gallot} for further results on
manifolds. Some isoperimetric estimates are also given in \cite{BK}.
The Cheeger constant has physical applications in modeling of land
slides \cite{CCCNP}, in image processing \cite{Ionuscu} and in
fracture mechanics \cite{Keler}.

 Any set $F$ in $\O$ which realizes the infimum
\eqref{eq: Cheeger} is called a Cheeger set in $\O$. If $\O$ is a
minimizer for \eqref{eq: Cheeger}, we  say that $\O$ is \emph{self-Cheeger}. $\O$ will be called
 \emph{uniquely self-Cheeger}, if  any  Cheeger  set in $\O$  is equal to $\O$ up to a  Riemann measure zero set.
 The main result of this paper is the following.
\begin{Theorem}\label{theo0}
Let $(\M,g)$ be a smooth compact Riemannian manifold  of dimension
$N\geq2$. Then there exists a family of uniquely self-Cheeger sets
$(\O_\e)_{\e\in(0,\e_0)}$ with
 $$h(\O_\e)=\frac{N}{\e}\quad\textrm{for all}\quad \e\in(0,\e_0).$$
Moreover, if $p_0$ is a non-degenerate critical point of the scalar
curvature  $s$ of $g$, then the family $(\de \O_\e)_{\e\in
(0,\e_0)}$ constitutes a smooth foliation of a neighborhood of
$p_0$.
\end{Theorem}
Theorem \ref{theo0} relies  on  Theorem \ref{the3} due
to M.M. Fall and the author in \cite{FallIgnace} together with
Theorem \ref{theo2} below.  Indeed Theorem \ref{theo2} gives a
sufficient condition on the Ricci curvature for a domain
$\O\subset\M$  to be a uniquely self-Cheeger while Theorem
\ref{the3} provides a family of domains where this condition is
always satisfied.
\begin{Theorem}\label{theo2}
Let $( \M,g)$ be a smooth compact Riemannian manifold of dimension
$N\geq2$ and $\O$ a bounded  domain  in $\M$ with Lipschitz
boundary. Assume that there exists a function $u\in C^2(\ov{\O})$
such that
\begin{equation}\label{eq:Pro}
  \begin{cases}
    -\Delta_g u=1& \quad \textrm{ in } \quad\Omega\vspace{3mm}\\
u=0&  \quad\textrm{ on }\quad \partial\Omega\vspace{3mm}\\
g(\n u , \eta)=-\l&   \quad \textrm{ on }\quad\partial\Omega,
  \end{cases}
  \end{equation}
where $\l$ is a positive constant,
$\Delta_{g}=\textrm{div}(\nabla_{g})$ is the Laplace-Beltrami
operator and  $\eta$ is the unit outer  normal of $\partial\Omega$.
If
\begin{equation}\label{eq: assum Ricci}
{Ric}_g(\nabla u,\nabla u)>-\frac{1}{N}\quad
\textrm{in}\quad \O,
\end{equation}
 then $\Omega$ is a uniquely self-Cheeger set
with Cheeger constant equal to $\dfrac{1}{\l}$.
\end{Theorem}
\begin{Theorem}\label{the3}
 Let $(\M,g)$ be a compact Riemannian manifold of dimension $N\geq2$. Then there exists $\e_{0}>0$ such that
 for all $\e\in (0,\e_0)$,  there
exist a smooth domain $\Omega_{\e}$ and a  function $u_{\e}\in
C^2(\ov{\O_\e})$  such that
\begin{align}\label{eq:problemz}
  \begin{cases}
    -\Delta_{g} u_\e=1& \quad \textrm{ in } \quad \Omega_\e \vspace{3mm}\\
u_\e=0&  \quad\textrm{ on }\quad\partial\Omega_\e \vspace{3mm}\\
{g}(\nabla_{ {g}} {u_\e}, {\nu}_\e)=-\dfrac{\e}{N}&   \quad \textrm{
on }\quad\partial\Omega_\e.
  \end{cases}
  \end{align}
  In particular if $p_0$ is a non-degenerate critical point of the scalar
curvature $s$ of $g$, then the  family $(\de \O_\e)_{\e\in
(0,\e_0)}$ constitutes a smooth foliation of a neighborhood of
$p_0$.
  \end{Theorem}
The domains $\Omega_\e$ in Theorem \ref{the3} are perturbations of
geodesic balls  $B_\e^g(p)$ centered at $p$ with radius $\e$.
Indeed, we have  that $ \Omega_{\e} =(1+v^{p,\e}) B_\e^g(p)$, with
$v^{p,\e}:\de B_\e^g(p) \to \R $ satisfying
\begin{equation}\label{eq: est}
||v^{p,\e}||_{C^{2,\alpha}(\partial B^{g}_{\e}(p))}\leq c\e^{2},
\end{equation}
where the constant $c$ is independent on $\e$.

If $p_0$ is a non-degenerate critical point of the scalar curvature
$s$ of $g$, we obtain a precise form of the boundary of the domains
$\Omega_\e$: there exists a smooth function $\o^{\e}: S^{N-1}\to
\R_+$ such hat
$$
\de \O_\e=\left\{\exp_{p_0}\left(\o^{\e}(y)\sum_{i=1}^N y^i
E_i\right),\quad y\in S^{N-1} \right\}
$$
and moreover the map $\e\mapsto \o^\e$ satisfies
$\de_\e\o^\e|_{\e=0}=1$.\\

From the result of Micheletti and Pistoia in \cite{MiPi}, it is
known that for a generic metric on a manifold, all critical points
of the scalar curvature are non-degenerate. The result of Theorem
\ref{theo0} then implies that for a generic metric $g'$, a
neighborhood of  any critical point of the scalar curvature can be
foliated by a  family $(\de \O_\e)_{\e\in(0,\e_0)}$, where $\O_\e$
is a uniquely self-Cheeger with Cheeger constant equal to $N/\e$.
This can be seen as the analogue  of Ye's result, who proved in
\cite{Ye} that around a non-degenerate critical point of the scalar
curvature function in a Riemannian manifold, there exists a unique
regular foliation by constant mean curvature spheres. It is an open
question to know whether the hypersurfaces constructed by Ye enclose
self-Cheeger sets. Indeed, Nardulli \cite{SN} proved that the sets
constructed by Ye are solutions to the isoperimetric problem in a
small neighborhood of the non-degenerate maxima of the scalar
curvature function. From the author's knowledge, except in space of
constant sectional curvature, a condition under which a solution to
the isoperimetric problem is a self-Cheeger set is not yet known.
The relation  between solutions to the isoperimetric problem and
self-Cheeger sets is discussed in Remark \ref{eq: rem}.

\section{Proof of Theorem  \ref{theo2}}
Let $u$ be the solution of \eqref{eq:Pro} in $\O$. We claim that \be
\label{eq:est-nab-us}
|\nabla u|_g< \l \quad \textrm{ in \quad
$\O$.} \ee

Indeed, we consider the smooth function $P: \M\rightarrow
\mathbb{R}$ defined  by
$$
P(q):=|\nabla u(q)|^2,\quad q\in \M.
$$
We recall Bochner's formula (see \cite{GSa}): for all $u\in
C^3(\M)$,
\begin{equation}\label{eq: Bocner}
\Delta_g|\nabla u|_g^2=2|\textrm{Hess}_g(u)|^2+ 2\langle \nabla u,
\nabla (\Delta_g u)\rangle_g+ 2\textrm{Ric}_g(\nabla u, \nabla u),
 \end{equation}
where $\textrm{Ric}_g$ is the Ricci curvature of $g$,
$\textrm{Hess}_g(u)$ stands for  the Hessian of $u$ and we have
denoted $|A|^2 = \textrm{tr}(AA^{t} )$ for a bilinear form $A$. By
the Cauchy-Schwarz inequality,
\begin{equation}\label{eq:deltea}
|\textrm{Hess}_g(u)|^2\geq\frac{1}{N}(\Delta_g u)^2.
\end{equation}
Since $-\Delta_{g}u=1$ in $\O$, we get from \eqref{eq: Bocner} and
\eqref{eq:deltea} that
$$
\frac{1}{2}\Delta_g  P(q) \geq\frac{1}{N}+\textrm{Ric}_g(\nabla
u(q), \nabla u(q)),\quad q \in \O.
$$
Using  the assumption $\textrm{Ric}_g(\nabla u(q),\nabla u(q))>-1/N$
in $\O$, we get $\Delta_g P>0\textrm{ in}\quad  \O$. Thus by the
strong maximum principle, $P$ attains  its maximum only on $\de \O$.
Since $P\equiv \l^2$ on $\partial \O$ by \eqref{eq:Pro}, the claim \eqref{eq:est-nab-us}
follows.

 Let $A\subseteq \O$ be  a set of finite perimeter.  Then
integrating \eqref{eq:Pro} on $A$ and  applying the claim \eqref{eq:est-nab-us}
with the Cauchy-Schwarz inequality, we get
\begin{equation}\label{eq cheeA}
|A|=-\int_{\de^{*} A} \langle \n u,\nu_A \rangle \textrm{d}\s_g \leq \int_{\de^{*} A} |\n u|\textrm{d}\s_g \leq\l
P( A).
\end{equation}
Recalling $\O$ is Lipschitz,  we integrate \eqref{eq:Pro} on  $\O$
 and use \eqref{eq cheeA} to get the inequality $$\frac{P(\O)}{|\O|}=\frac{1}{\l}\leq\frac{P(A)}{|A|}$$  which shows that  $
\O$ is self-Cheeger with Cheeger constant equal to $1/\l$.\\

To complete the proof, its remains to show that $\O$ is uniquely
self-Cheeger. Let $C$ be a Cheeger set inside $\O$,  we need to show
that $|\O\setminus C|=0$.

Since $C$ is a Cheeger set we have $|C|=\lambda P(C)$, and writing \eqref{eq cheeA} with $C$ in place of $A$, we get
$$
\int_{\de^{*} C} (\lambda-|\n u|)\textrm{d}\s_g=0.
$$
From  this and \eqref{eq:est-nab-us} we  get
$$|\n u|=\lambda \quad \textrm{on}\quad \partial^{*} C, \quad \sigma_g-\textrm{almost  everywhere}$$ and
hence $\partial ^{*} C\subseteq \partial \O$.

We now prove  that $|\O\setminus C|=0$. This follows from
Poincaré-Wirtinger's inequality \cite[Theorem 2.10]{Hebey}:
\begin{equation}\label{eq: poiWIR}
\forall w\in W^{1,1}(\O), \quad ||w-m(w)||_{L^1(\O)}\leq
k\int_{\O}|\n w|\textrm{d}v_g,
\end{equation}
 where
$k>0$ is a constant  only depending on $\O$ and
$$m(w)=\frac{1}{|\O|}\int_{\O} w \textrm{d}v_g.$$ Since $\O$ has Lipschitz boundary, for  every $w\in BV(\O)$ there exists a sequence $(w_n)\in
W^{1,1}(\O)\cap C^{\infty}(\O)$ such that
$$w_n\longrightarrow w \quad \textrm{in }\quad  L^1(\O)\quad \textrm{and}\quad \int _{\O} |\n w_n|\textrm{d}v_g \longrightarrow |Dw|(\O).$$
Inequality \eqref{eq: poiWIR} then also holds in the space $BV(\O)$.

For $w=1_C\in BV(\O)$, we have
$$
|1_{C}-m(1_{C})|=\frac{|\O \setminus C|}{|\O|}\quad \textrm{ in}
\quad C \quad \textrm{and} \quad
|1_{C}-m(1_{C})|=\frac{|C|}{|\O|}\quad \textrm{ in}\quad  \O
\setminus C.
$$
From this and \eqref{eq: poiWIR}, we get
$$
||1_{C}-m(1_{C})||_{L^1(\O)}=2\frac{|C|}{|\O|}|\O\setminus C| \leq k
P(C, \O).
$$
Since $\partial^* C\cap \O=\emptyset$, we have $P(C,\O)=|\partial^*
C\cap \O|=0$. Thanks to the above inequality, $|\O\setminus C|=0$,
and the proof of Theorem  \ref{theo2} is complete.\QED
\section{Proof of Theorem \ref{theo0}}
In view of Theorem \ref{theo2}, the proof of  Theorem \ref{theo0}
will be fulfilled once we show that the solution $u_\e$ to
\eqref{eq:problemz} satisfies
$$
\textrm{Ric}_g(\nabla u_\e, \nabla u_\e)>-\frac{1}{N}\quad
\textrm{in}\quad \O_\e.
$$
 Denote by
$\textrm{r}_{p}(\cdot):=\textrm{dist}_{g}(p,\cdot)$  the Riemannian
geodesic distance  function to $p$. We recall that for all
$\e\in(0,\e_0)$, the solution $u_\e$ to \eqref{eq:problemz} is given
(see  \cite[Section 3]{FallIgnace}) by
\begin{equation}\label{eq:psii}
u_\e(q)=\frac{\e^2-\textrm{r}_{p}^2(q)}{2N}+\e^4w_\e(q)\quad
\textrm{for all}\quad q\in \Omega_{\e},
\end{equation}
where
$$
|| w_\e||_{C^{2}(\M)}\leq c
$$
and $c$ is a positive  constant independent on $\e$.  Since $\M$ is
compact, there exists a positive constant $C$ depending only on $\M$
such that
\begin{equation}\label{eq Inefinadl}
\textrm{Ric}_g(\nabla u_\e(q), \nabla u_\e(q))\geq-C|\nabla
u_\e(q)|_g^2.
\end{equation}
Using \eqref{eq:psii}  and the fact that  $|\nabla
\textrm{r}_{p}(\cdot)|_g=1$, we get
$$
|\nabla u_\e(q)|_g^2\leq c\e^2,
$$
where $c$ is a positive constant only depending on $\M$. This
together with  \eqref{eq Inefinadl} implies that
\begin{equation}\label{eq Inefinadbl}
\textrm{Ric}_g(\nabla u_\e(q), \nabla u_\e(q))\geq-C\e^2,
\end{equation}
where $C$ is a positive constant only depending on $\M$. Finally, it
follows from \eqref{eq Inefinadbl} that for all $\e$ small enough
$$\textrm{Ric}_g(\nabla u_\e, \nabla
u_\e)>-\frac{1}{N}\quad \textrm{in $\O_\e$}.$$ This ends the proof
of Theorem  \ref {theo0}.\QED

\begin{Remark}\label{eq: rem}

It is an open problem to know if the sets enclosed by Ye's
hypersurfaces are self-Cheeger sets on $\cM$.  It is known (see
\cite{Druet-FK,SN}) that solutions  to the isoperimetric problem
with small  volume are near the maxima of the scalar curvature, and
if this maximum is non-degenarate, then
\cite{SN} implies that, solutions to the isoperimetric problem are among Ye's
sets.

As we shall see below,  it might be possible to have  some answers
to the above open problem, provided higher order Taylor expansions
of  the perimeter of an isoperimetric set, as volume tends to zero,
is available.

 In what follows, $B$ denotes the unit Euclidean ball of
$\mathbb{R}^N$ and we put $\o_{N-1}:=|\partial B|$.

The isoperimetric profile of a compact Riemannian manifold is given,
see \cite{Druet-FK,SN} by
$$
\mathcal{I}_{\M}(v,g)=\inf_{E\subset\M,~ |E|=v}P(E)
$$
while the Cheeger isoperimetric  profile is defined \cite[section
4]{mmf1}  by
\begin{equation}\label{eq:isop}
\mathcal{H}_{\M}(v,g)=\inf_{\O\subset\M,~ |\O|=v}h(\O),
\end{equation}
where $h(\O)$ is the Cheeger constant of the set $\O\subset\M$.

The profile $\mathcal{H}_{\M}$ was derived in \cite[Theorem
4.1]{mmf1} near $v=0$, ($v>0$)  and expands as
\begin{equation}\label{eq:chegerpro}
\mathcal{H}_{\M}(v,g)=\biggl(1+\mathcal{O}(v^{\frac{3}{N}})\biggl)\frac{\mathcal{I}_{\M}(v,g)}{v}
\end{equation}
and the isoperimetric profile of the manifold $\M$  has the
following expansion, \cite[Theorem 8]{SN}(assuming the absolute
maxima of the scalar curvature $s$ are non-degenerate critical
points),
$$
\mathcal{I}_{\M}(v,g)=c_Nv^{\frac{N-1}{N}}\biggl(1-\frac{S}{2N(N+2)}\biggl(\frac{v}{|B|}\biggl)^{\frac{2}{N}}+\mathcal{O}(v^{\frac{4}{N}})\biggl),
$$
where
$$S=\max_{p\in\M}s(p)\quad\textrm{and}\quad
c_N=\frac{\o_{N-1}}{|B|^{\frac{N-1}{N}}}.$$

Let  $E$  be a solution to the isoperimetric problem with volume
$v$, then by \eqref{eq:chegerpro} we get
\begin{equation}\label{eq:chegerprddo}
\mathcal{H}_{\M}(v,g)=\biggl(1+\mathcal{O}(v^{\frac{3}{N}})\biggl)\frac{P(E)}{|E|}.
\end{equation}
 If we further
assume  that $E$ is a self-Cheeger set so that
$$h(E)=\frac{P(E)}{|E|},$$ then by \eqref{eq:chegerprddo} and the
definition of $\mathcal{H}_{\M}$ in \eqref{eq:isop}, we get
$$
\biggl(1+\mathcal{O}(v^{\frac{3}{N}})\biggl)h(E)=\mathcal{H}_{\M}(v,g)\leq
h(E).$$ A higher order expansion of the perimeter of an
isoperimetric set is then needed in order to derive a (necessary)
condition  under which a solution to the isoperimetric problem (or
more generally the enclosure of a Ye's hypersurface) is a
self-Cheeger set.
\end{Remark}

\bigskip
\noindent \textbf{ Acknowledgements}:  The author thanks his
supervisor Mouhamed Moustapha Fall for his helpful suggestions and
insights throughout the writing of this paper. He is also grateful
to the reviewer for his valuable comments.  Part of this work was
done when the author  was visiting the Institute of Mathematics of
the Goethe-University Frankfurt. This work is supported by the
German Academic Exchange Service (DAAD).

\end{document}